\documentclass[12pt]{article}

\usepackage{amsmath}
\usepackage{amssymb}
\usepackage{amscd}
\newcommand{\qed}{\hbox{\unskip\nobreak\hfil
        \penalty50\hskip1em\hbox{}\nobreak\hfil
        $\square$\parfillskip=0pt\finalhyphendemerits=0 \par}}

\newtheorem{dfn}{Definition}[section]
\newtheorem{rem}[dfn]{Remark}

\newtheorem{thm}[dfn]{Theorem}
\newtheorem{lem}[dfn]{Lemma}
\newtheorem{prop}[dfn]{Proposition}
\newtheorem{cor}[dfn]{Corollary}

\newtheorem{defn}[dfn]{Definition}

\def\no{\noindent}

\def\={=\joinrel =}
\def\0{\emptyset}

\def\proof{\par\medskip\noindent{\it Proof: }}
\def\lra{\longrightarrow}

\def\>{\rangle}
\def\<{\langle}

\def\C{\mathbb C}

\def\t{\tilde}
\def\R{\mathbb R}
\def\P{\mathbb P}
\def\Z{\mathbb Z}
\def\E{{\mathbb E}}
\def\p{{\cal P}}

\def\Q{\mathbb Q}
\def\F{{\mathbb F}}

\def\M{{\cal M}}
\def\eps{\epsilon}
\def\H{\mbox{\rm H}}

\def\al{\alpha}

\def\ga{\gamma}
\def\Ga{\Gamma}
\def\del{\delta}

\def\si{\sigma}
\def\L{{\cal L}}

\def\la{\lambda}

\def\8{\infty}
\def\Om{\Omega}
\def\om{\omega}
\def\D{\partial}

\def\e{{\cal E}}

\def\hook{\hookrightarrow}

\def\1{\hbox{\bf 1}}
\def\m{{\mathfrak m}}

\def\BE{\begin{equation}}
\def\EE{\end{equation}}
\def\bul{\bullet}

\def\Hom{\mbox{\rm Hom}}

\def\p{{\mathcal P}}
\begin{document}
\title{Quantization of bending deformations of polygons in $\E^3$, 
hypergeometric integrals and the Gassner representation }
\author{Michael Kapovich \ 
and John J.\ Millson}
\date{February 17, 2000}
\maketitle

\begin{abstract}
\no The Hamiltonian potentials of the bending deformations of 
$n$-gons in $\E^3$ studied in \cite{KM} and \cite{Kl} give rise to a 
Hamiltonian action of the Malcev Lie algebra $\p_n$  of the pure braid group 
$P_n$ on the moduli space $M_r$ of $n$-gon linkages with the side-lengths 
$r= (r_1,..., r_n)$ in $\E^3$. If $e\in M_r$ is a singular point we may 
linearize the vector fields in $\p_n$ at $e$. This linearization yields a flat connection 
$\nabla$ on the space $\C^n_*$ of $n$ distinct points on $\C$. We show 
that the monodromy of $\nabla$ is the dual of a quotient of a specialized reduced Gassner representation.  

\end{abstract}

\section{Introduction} 

In \cite{KM} and \cite{Kl} certain Hamiltonian flows on the moduli space 
$M_r$ of $n$-gon linkages in $\E^3$ were studied. In \cite{KM} these flows 
were interpreted geometrically and called {\em bending deformations of 
polygons}. In \cite{Kl}, Klyachko pointed out that the Hamiltonian 
potentials of the bending deformations gave rise to a Hamiltonian action 
of $\p_n$, the Malcev Lie algebra of the pure braid group $P_n$ 
(see \S 3), on $M_r$. 
It is a remarkable fact, see \cite[Lemma 1.1.4]{K1}, that a representation 
$\rho: \p_n\to End(V)$, $dim(V)<\infty$, gives   rise to a flat connection 
$\nabla$ on the vector 
bundle $\C^n_*\times V$ over $\C^n_*$, the space of distinct points in 
$\C$. Accordingly the monodromy representation of $\nabla$ yields a 
representation $\hat\rho: P_n\to Aut(V)$. 

We see then that if we can find a finite dimensional representation of 
the Lie algebra ${\cal B}\subset C^{\infty}(M_r)$ generated by the 
bending Hamiltonians under the Poisson bracket , i.e. if we can 
``quantize''   ${\cal B}$, then we 
will obtain a representation of $P_n$. Klyachko suggested using a geometric 
quantization of $M_r$ to quantize ${\cal B}$. This appears to be difficult to 
carry out because the bending flows do not preserve a polarization. Note 
however that the problem of quantizing a Poisson subalgebra of $C^{\infty}(M_r)$ 
can be solved immediately if the functions in the subalgebra have a 
common critical point 
$x\in M_r$. For in this case we may simultaneously linearize all the 
Hamiltonian fields at $x$. We are fortunate that simultaneous critical 
points  for the algebra ${\cal B}$ exist if $M_r$ is singular. 
Indeed, a degenerate $n$-gon (i.e. an $n$-gon which is 
contained in a line $L$) is a critical point of all bending Hamiltonians. 

The point of this paper is to compute the representation 
$\hat\rho_{\eps,r}: P_n\to Aut(T_{\eps,r})$ associated to a 
degenerate $n$-gon $P$. Here 
$T_{\eps,r}= T_{P}(M_r)$ and 
$\eps= (\eps_1,..., \eps_n)$, $\eps_i\in \{\pm 1\}$, and 
$r=(r_1,..., r_n), r_i\in \R_+$, are defined as follows. Fix an orientation 
on $L$. The number $r_i$ is the length of the $i$-th edge of $P$. 
Define $\eps_i$ to be $+1$ if the $i$-th edge is positively oriented and 
$\eps_i=-1$ otherwise. We call $\eps  =(\eps_1,..., \eps_n)$   
the vector of edge-orientation of $P$. 

Our formula for $\rho_{\eps,r}: \p_n \to End(T_{\eps,r})$ is in terms 
of certain $n\times n$ 
matrices $J_{ij}(\la)$ which are called {\em Jordan-Pochhammer} 
matrices. Let $\la=(\la_1,...,\la_n)$ be an $n$-tuple of 
complex numbers. 
Define matrices 
  $J_{ij}(\la)$ for $1\le i < j\le n$ by

$$
 \bordermatrix{&~&i-\hbox{\rm th~column}&~&j-\hbox{\rm th~column}&~\cr
~&0\ldots 0 &0&0\ldots 0&0&0\ldots 0\cr 
i-\hbox{\rm th~row}&0\ldots 0 &\la_j &0\ldots 0 &-\la_j&0\ldots 0\cr
~&0\ldots 0 &0&0\ldots 0&0&0\ldots 0\cr 
~&\vdots&\vdots&\vdots&\vdots&\vdots\cr 
j-\hbox{\rm th~row}&0\ldots 0 &-\la_i &0\ldots 0 &\la_i&0\ldots 0\cr
~&0\ldots 0 &0&0\ldots 0&0&0\ldots 0\cr 
~&\vdots&\vdots&\vdots&\vdots&\vdots\cr}= J_{ij}(\la).
$$

Define $J_{ii}=0$ and $J_{ij}(\la)= J_{ji}(\la)$ for $i>j$. We have 
(as can be verified easily)

\begin{lem}
The matrices  $\{J_{ij}(\la)\}$ satisfy the infinitesimal braid relations:
\begin{itemize}
\item $[J_{ij}(\la), J_{kl}(\la)]=0$ if $\{i,j\}\cap \{k,l\}=\emptyset$. 
\item $[J_{ij}(\la), J_{ij}(\la) + J_{jk}(\la) + J_{ki}(\la)]=0$, $i,j,k$ are 
distinct. 
\end{itemize}
\end{lem}

Consequently the assignment $\rho_{\la}(X_{ij})= J_{ij}(\la)$ (see Section 
3 for the meaning of $X_{ij}$) yields a representation 
$\rho_{\la}: \p_n\to M_n(\C)$ and a flat connection $\nabla$ on 
$\C^n_* \times \C^n$. Here we realize $\C^n$ as the space of {\em row} vectors with $n$ components. 
It is immediate that the subspace 
$\C^n_0\subset \C^n$ defined by 
$$
\C^n_0=\{z\in \C^n : \sum_i z_i=0\}
$$
is invariant under $\rho_\la$, in fact $\rho_\la(\p_n)(\C^n)\subset \C^n_0$. 
Now we assume $\sum_{i=1}^n \la_i=0$. Then $\la\in \C^n_0$ and we see that 
$\rho_{\la}(\p_n)(\la)=0$. Thus we have a $\p_n$-invariant filtration 
$$
\C \la \subset \C^n_0 \subset \C^n .
$$
Define $W_\la=\C^n_0/\C\la$. Now let $P$ be a degenerate $n$-gon with 
side-lengths 
$r=(r_1,...,r_n)$ and edge-orientations $\eps= (\eps_1,...,\eps_n)$. 
Our first main theorem is

\medskip
{\bf Theorem A.}  {\em There is a $\p_n$-invariant almost complex structure 
$J^{\eps}$ on $T_{\eps,r}$ such that there is an isomorphism of 
$\p_n$-modules 
$T^{1,0}_{\eps,r}\cong W_\la$ for $\la:= (\sqrt{-1}\eps_1 r_1,..., 
\sqrt{-1}\eps_n r_n)$.}  

\medskip
Here $T^{1,0}_{\eps,r}= \{w\in T_{\eps,r}\otimes \C: J^{\eps}w= \sqrt{-1}w\}$. 
We have

\medskip
{\bf Corollary.} {\em The flat connection on $\C^n_* \times  T^{1,0}_{\eps,r}$ 
has the connection form
$$
\om = \sum_{1\le i< j\le n} \frac{dz_i -dz_j}{z_i-z_j} \otimes  J_{ij}(\la)
$$
with $\la$ as above.}

\medskip
We then adapt the methods of \cite{K1} to give formulae for multivalued 
parallel sections of $\nabla$ in terms 
of hypergeometric integrals and to compute the monodromy of $\nabla$. 

Before stating our first formula for the monodromy of $\nabla$ we need more 
notation. Let $\ga_j$, $1\le j\le n$, be the free generators of the free group $\F_n$. Define 
the character $\chi: \F_n \to \C^*$ by $\chi(\ga_j)= e^{2\pi i \la_j}$, $1\le j\le n$ 
(recall that $\la_j= \sqrt{-1}\eps_j r_j$). Let $\C_{\chi^{-1}}$ be the 1-dimensional module 
(over $\C$) in which the free group $\F_n$ acts by $\chi^{-1}$. The pure braid 
group $P_n$ acts by automorphisms on $\F_n$ so that the character $\chi$ is fixed. 
Thus we have the associated action of $P_n$ on $\H_1(\F_n,  \C_{\chi^{-1}})$. 
We let $\Ga_n= \pi_1(\C\P^1 - \{z_1,...,z_n\})$ be the fundamental group 
of the $n$ times punctured sphere. Hence $\Ga_n$ is the quotient of $\F_n$ by 
the normal subgroup generated by $\ga_1\ldots \ga_n$. Since 
$\chi(\ga_1 \ldots \ga_n)=1$, the character $\chi$ induces a character of $\Ga_n$. 
The group $P_n$ fixes $\ga_1\ldots \ga_n$ and consequently acts on $\Ga_n$ and on 
$\H_1(\Ga_n, \C_{\chi^{-1}})$. We can now state

\medskip
{\bf Theorem B.} {\em The monodromy representation of 
$\nabla$ is equivalent to the 
representation of $P_n$ on $\H_1(\Ga_n,  \C_{\chi^{-1}})$.}

\medskip 
In \S \ref{Gassner} we define the {\em Gassner representation} of the pure braid group, 
the {\em reduced Gassner representation} and their {\em specializations} 
via characters of the free group. Let $\L$ is the $\C$-algebra of Laurent 
polynomials on $t_1,..,t_n$.

\medskip
{\bf Theorem C.}  {\em The monodromy representation of $\nabla$ is dual 
to the quotient of the reduced Gassner representation 
 $Z^1(\Ga_n,\L)$  
specialized at $t_j= e^{-2\pi \eps_j r_j}$, where we quotient by the  
1-dimensional subspace $B^1(\Ga_n, \C_{\chi})$ fixed by $P_n$.} 

\medskip
Our results appear to be related to those of \cite{DM} and \cite{L} but there are significant 
differences. In \cite{L}, D.~D.~Long linearizes the action of $P_n$ 
on the moduli space of $n$-gon linkages in $S^3$ obtained from the action of $P_n$ on 
$$
\Hom(\pi_1(S^2-\{z_1,...,z_n\}), SU(2))/SU(2)
$$
by precomposition. The corresponding action of $P_n$ on $M_r$ 
is trivial in our case, see \cite[Remark 5.1]{KM}. In \cite{DM}, Deligne and 
Mostow arrive at the Gassner representation by considering a variation of Hodge 
structure over $\C^n_*/PGL_2(\C)\subset M_r$. They obtain the quotient 
(by the 1-coboundaries) of the 
reduced Gassner representation specialized at 
$(e^{2\pi i r_1},...,e^{2\pi i r_n})$; 
we obtain the dual of the  quotient of the reduced  
Gassner representation specialized at 
$(e^{-2\pi \eps_1 r_1},...,e^{-2\pi \eps_n r_n})$. 
Here we must assume $\sum_{i=1}^n r_i=2$ to be consistent with \cite{DM}. 
Our representation lies in $GL(n-2, \R)$; their 
representation is in $U(n-3,1)$.

\medskip
{\bf Acknowledgements.} It is our pleasure to thank Ragnar 
Buchweitz and Richard Hain for helpful conversations.

\section{The moduli space of $n$-gon linkages in $\E^3$.} 

Let $Pol_n(\E^3)$ be the space of (closed) $n$-gons with distinguished vertices 
in the Euclidean space $\E^3$. An $n$-gon $P$ is defined to be an ordered 
$n$-tuple of points $(v_1,..., v_n)\in (\E^3)^n$. The point $v_i$ is called 
the $i$-th vertex of $P$. The vertices are joined in cyclic order by edges 
$e_1,..., e_n$ where $e_i$ is the oriented segment from $v_i$ to $v_{i+1}$.  
We think of $e_i$ as a vector in $\R^3$. Two polygons $P= (v_1,..., v_n)$ 
and $Q= (w_1,..., w_n)$ are identified if and only if there exists an orientation-preserving isometry 
$g$ of $\E^3$ such that $g(v_i)=w_i$, $1\le i\le n$. Let $r=(r_1,...,r_n)$ 
be an $n$-tuple of positive real numbers. Then $M_r$ is defined to be the 
moduli space of $n$-gons with the side-lengths $r_1,..., r_n$ modulo 
isometries as above. An element of $M_r$ will be called a 
{\em closed $n$-gon linkage}. 

We will also need the moduli space space $N_r$ of ``open'' 
$n$-gon linkages. To obtain $N_r$ we repeat the above construction of $M_r$ 
except we do not assume the end vertex $v_{n+1}$ of the edge $e_n$ is equal to $v_1$.

The starting point of \cite{KM} was the observation that 
$$
M_r= \{e= (e_1,...,e_n) \in \prod_{i=1}^n S^2(r_i) : 
e_1 + ... + e_n =0\}/SO(3).
$$  
This equality exhibits $M_r$ as the symplectic quotient of 
$\prod_{i=1}^n S^2(r_i)$ and has many consequences. First $M_r$ is a 
complex analytic space with isolated (quadratic) singularities. The 
smooth part of $M_r$ 
is a K\"ahler manifold. The singular points of $M_r$ are the equivalence 
classes 
of degenerate $n$-gons. Thus $M_r$ is singular if and only if $r$ is the 
set of side-lengths of a degenerate $n$-gon.   
 
In \cite{KM} we introduced bending deformations of closed polygonal 
linkages in $\E^3$, see also \cite{Kl}. Suppose $P= e=(e_1,..., e_n)$. Let 
$I\subset \{1,...,n\}$ be a subset and define 
$f_I\in C^{\infty}(M_r)$ by 
$$
f_I(e)= \|\sum_{i\in I} e_i\|^2.
$$
Then $f_I$ is the Hamiltonian potential of a Hamiltonian vector field $B_I$. 
The vector $e_I= \sum_{i\in I} e_i$ is constant along an integral curve of 
$B_I$. By \cite[Lemma 3.5]{KM},  $B_I(e)= (\del_1,..., \del_n)$, where 
$\del_i= e_I\times e_i$, $i\in I,$ and $\del_i=0$ for $i\notin I$.  
The integral curves of $B_I$ are obtained as follows. Define an element 
$ad(e_I) \in so(3)$ by
$$
ad(e_I)(v)= e_I \times v
$$ 
and a one-parameter group $R_I(t)\subset SO(3)$ by
$$
R_I(t)= \exp (t \ ad (e_I)).
$$
Then the integral curve $e(t)$ of $B_I$ passing through $e$ is given by
$$
e_i(t)= R_I(t) e_i, i\in I
$$
$$
e_j(t)= e_j, j\notin I.
$$
This motion of a polygon $P$ has a simple geometric interpretation if the 
elements of $I$ are consecutive. In this case $e_I$ is a diagonal and it 
divides the polygon into two parts. Keep one part fixed and {\em bend} the 
polygon by rotating the other part around the diagonal with the angular 
speed  $\|e_I\|$. For this reason we call the above motion a bending 
deformation of the polygon. We will be specifically interested in the case 
$I=\{i,j\}$, $i< j$. We abbreviate $f_{\{i,j\}}$ to $f_{ij}$ 
and $B_{\{i,j\}}$ to $B_{ij}$ . We have:
$$
f_{ij}(e)= \|e_i + e_j\|^2.
$$ 

\begin{lem}
Let $e\in M_r$ be a degenerate polygon. Then $B_{ij}(e)=0$ for all $i, j$. 
\end{lem}
\proof The bending field $B_{ij}$ is given by
$$
B_{ij}(e)= (0,..., (e_i + e_j)\times e_i, 0,..., 
(e_i + e_j)\times e_j, 0,...)= 
(0,..., e_j\times e_i, 0,..., e_i \times e_j, 0,...).
$$
If $e$ is degenerate then $e_i$ and $e_j$ are linearly dependent, so 
$e_i \times e_j=0$.  \qed

\begin{rem}
In fact $B_I(e)=0$ for all $I$ if $e$ is degenerate. 
\end{rem}

Define $\t N_r:= \prod_{i=1}^n S^2(r_i)$ where $S^2(r_i)$ is the round 
2-sphere of the radius $r_i$. We also define $\t M_r\subset  \t N_r$ by 
$$
\t M_r =\{ e\in \t N_r : \sum_{i=1}^n e_i =0\}.
$$
Hence $N_r$ is the quotient of $\t N_r$ by $SO(3)$ and $M_r$ is the 
quotient of $\t M_r$ by $SO(3)$.

\section{The Malcev Lie algebra of the pure braid group.}

Let $P_n$ be the pure braid group on $n$ strands in $\C$ 
(see \cite[\S 1]{C}). Let $\C^n_*$ denote the subset of $\C^n$ 
consisting of distinct $n$-tuples. Then $P_n$ 
is isomorphic to the fundamental group of $\C^n_*$.  

Let $\p_n$ be the Malcev Lie algebra of $P_n$, see \cite{ABC}. Kohno found 
the following presentation for $\p_n$ in \cite{K2} 
(see also \cite[Proposition 3.2.1]{I}). 

\begin{lem}
\label{3.1}
The Lie algebra $\p_n$ is the quotient of the free Lie algebra over $\Q$ 
generated by $X_{ij}, 1\le i , j\le n$, subject to the relations:
\begin{enumerate}
\item $X_{ii}=0$, $1\le i\le n$. 
\item $X_{ij}= X_{ji}$, $1\le i , j\le n$
\item $[X_{ij}, X_{kl}]=0$ if $\{i,j\}\cap \{k,l\}=\emptyset$. 
\item $[X_{ij}, X_{ij} + X_{jk} + X_{ki}]=0$, $i,j,k$ are distinct. 
\end{enumerate}
\end{lem}

We will now see that any finite dimensional representation of $\p_n$ 
induces a finite dimensional representation of $P_n$ on the same vector 
space. This remarkable fact is an immediate consequence of the following 
lemma of Kohno \cite[Lemma 1.1.4]{K1}.

\begin{lem}
Suppose $V$ is a finite dimensional vector space and $A_{ij}, 1\le i,j\le n$, 
are elements of $End(V)$ such that $A_{ii}=0$ and $A_{ij}=A_{ji}$.  Let 
$\nabla$ be the connection on the trivial $V$ bundle over $\C^n_*$ with 
connection form
$$
\om= \sum_{1\le i< j\le n} \frac{dz_i - dz_j}{z_i -z_j} \otimes A_{ij}.
$$ 
Then $\nabla$ is flat if and only if the relations (3) and (4) for $\p_n$ 
are satisfied by the $A_{ij}$'s. 
\end{lem}

Thus there is a 1-1 correspondence between Lie algebra homomorphisms 
$\rho: \p_n \to End(V)$ and flat connections $\nabla$ on $\C^n_*\times V$ 
of the above form. Suppose we are given $\rho$ as above. Since 
$\pi_1(\C^n_*, z)\cong P_n$ ($z$ is a base-point), the monodromy 
representation of $\nabla$ gives an induced representation of $P_n$ to $Aut(V)$. 

Let $F: \C^n_* \to V$ be a smooth map. Then $F$ induces a parallel 
section of $\nabla$ if and only if $F$ satisfies the equation (of the 
$V$-valued 1-forms on $\C^n_*$)
$$
dF= \sum_{1\le i< j\le n} \frac{dz_i - dz_j}{z_i -z_j} \otimes A_{ij}(F).
$$

\section{A Hamiltonian action of $\p_n$ on $M_r$.}

We define the function $f_{ij}$ on $\t{N}_r$ by  
$$
f_{ij}(e)= \|e_i + e_j\|^2.
$$
The next proposition was proved in \cite{Kl}. Since it is central to 
our paper we give a proof here.

\begin{prop}
\label{4.1}
\begin{enumerate}
\item $f_{ij}= f_{ji}$. 
\item $\{ f_{ij}, f_{kl}\}=0$, if $\{i,j\}\cap \{k,l\}=\emptyset$. 
\item $\{f_{ij}, f_{ij} + f_{jk} + f_{ki}\}=0$, if $i,j,k$ are distinct. 
\end{enumerate}
\end{prop}
\proof The assertions (1) and (2) are obvious. The third assertion will be a 
consequence of the following discussion. Since $\t{N}_r$ is a symplectic leaf 
of the Lie algebra $(\R^3, \times)$ equipped with the Lie Poisson structure it 
suffices to prove (3) for the functions $f_{ij}$ extended to $(\R^3)^n$ using 
the same formula. Let $g_{ij}: (\R^3)^n\to \R$ be given by $g_{ij}(e)= 
e_i \cdot e_j$ 
and $h_{ijk}: (\R^3)^n\to \R$ be given by $h_{ijk}(e)= e_i \cdot (e_j 
\times e_k)$. 

\begin{lem}
$\{g_{ij}, g_{jk}\}= -h_{ijk}$. 
\end{lem}
\proof It suffices to prove the lemma for $i=1, j=2, k=3$. We use coordinates 
$(x_i, y_i, z_i), 1\le i\le n$, on $(\R^3)^n$. Then
$$
\{ x_i, y_i\}= z_i, \{y_i, z_i\}= x_i, \{z_i, x_i\}= y_i, 1\le i\le n\quad .
$$
We have
$$
\{ g_{12}, g_{23}\} = \{x_1 x_2 + y_1 y_2 + z_1 z_2, x_2 x_3 + y_2y_3 + 
z_2 z_3\}= 
$$
$$
= \{ x_1 x_2, y_2 y_3\} + \{ x_1 x_2, z_2 z_3\} + \{y_1y_2, x_2 x_3\} + 
\{y_1y_2, z_2 z_3\} + \{z_1z_2, x_2 x_3\} +\{z_1z_2, y_2 y_3\}=$$
$$
=x_1 y_3 z_2 - x_1 y_2 z_3 - x_3 y_1 z_2 + x_2 y_1 z_3 + x_3 y_2 z_1 - 
x_2 y_3 z_1 = -e_1 \cdot (e_2 \times e_3).\quad \qed
$$

\begin{cor}
$\{f_{ij}, f_{jk}\} = -4 e_i \cdot (e_j \times e_k)$. 
\end{cor}
\proof $f_{ij}= f_{ii} + f_{jj} + 2g_{ij}$. But $f_{ii}$ and $f_{jj}$ are 
Casimirs. 
\qed

We now prove the 3-rd assertion. 
$$
\{f_{ij}, f_{ij} + f_{jk} + f_{ki}\} = \{ f_{ij}, f_{jk}\}  + 
\{ f_{ij}, f_{ki}\}= \{ f_{ij}, f_{jk}\} + \{ f_{ji}, f_{ik}\}=$$
$$
\{ f_{ij}, f_{jk}\} + \{ f_{ji}, f_{ik}\}= -4e_i \cdot e_j \times 
e_k - 4e_j \cdot e_i \times e_k = - 4e_i \cdot e_j \times e_k + 4e_i 
\cdot(e_j \times e_k)=0.
$$
\qed

Since the function $f_{ij}$ is $SO(3)$-invariant it induces a function 
(which is again denoted by $f_{ij}$) on $M_r$. The Poisson bracket of 
these functions remain the same and we obtain

\begin{thm}
There exists a Hamiltonian action of the Lie algebra $\p_n$ on the 
symplectic manifold $\t{N}_r$. This action induces an action on $M_r$. 
\end{thm}

>From Lemma \ref{3.1} and Proposition \ref{4.1} we see that if we can 
find a finite-dimensional representation of the Lie subalgebra of 
$C^{\infty}(M_r)$ generated by $\{f_{ij}, 1\le i < j\le n\}$ then we 
will get a representation of $P_n$. As explained in the introduction 
we obtain such a representation on $T_e(M_r)$ for a degenerate $n$-gon $e$. 

\section{Linearization of the bending fields at degenerate polygons.} 

This section is the heart of the paper. We compute $A_{ij}\in 
End(T_e(M_r))$, the linearization of the bending field $B_{ij}$ 
at a degenerate polygon $e\in M_r$. 
%By \cite[Lemma 3.5]{KM} the bending field $B_{ij}$, $1\le i<j\le n$, 
%is given by  $$
%B_{ij}(e)=  (0,..., e_j\times e_i, 0,..., e_i \times e_j, 0,...)$$
%where the nonzero entries are in the $i$-th and $j$-th spots. 
Now assume that $e$ is degenerate, so we may write
$$
e= (r_1 \eps_1 u,..., r_n\eps_n u)
$$
for some vector $u\in S^2$ and $\eps_i= \pm 1$.

Let $M$ be a manifold, $m\in M$. We recall the definition of the 
linearization $A_X \in End(T_m(M))$ of a vector field $X$ at a 
point $m$ where $X(m)=0$. Choose 
a connection $\nabla$ on $T(M)$. Let $u\in T_m(M)$, then
$$
A_X(u):= (\nabla_u X)(m)
$$
Since $X(m)=0$, $A_X$ is independent of the choice of connection. 

For the case in hand the above definition must be modified since $M_r$ 
is singular at $e$. There is a commutative algebra version of the above 
construction that goes as follows. Assume $M$ is a real affine variety, 
$m\in M$ and  $X$ is a vector field on $M$ satisfying $X(m)=0$. Let $\m$ 
be the maximal ideal of $m$. Then (since $X(m)=0$) we have $X\m \subset 
\m$ whence $X\m^2 \subset \m^2$ and $X$ induces an element of 
$End(\m/\m^2)= End(T_m^*(M))$. By duality we obtain $A_X\in End(T_m(M))$. 
The reader will verify that if $m$ is a smooth point of $M$ then the two 
definitions coincide. 

We now compute the linearization of $B_{ij}$ at $e$ in $M_r$. Recall that we have a diagram
$$
\begin{array}{ccc}
\t{M}_r & \longrightarrow & \t{N}_r\\
\downarrow & ~ & \downarrow\\
{M}_r & \longrightarrow & {N}_r
\end{array}
$$ 
where $\t{N}_r= S^2(r_1)\times ... \times S^2(r_n)$ and $\t{M}_r= \{e\in 
\t{N}_r : \sum_{i=1}^n e_i =0\}$. Define $g_{ij}: \t{N}_r \to \R$ by $g_{ij}(e)= \|e_i + e_j\|^2$. 
Hence $g_{ij}|\t{M}_r$ is $SO(3)$-invariant and descends to the function $f_{ij}$ on $M_r$. 
Let $\t{B}_{ij}$ be the Hamiltonian vector field of $g_{ij}$. Then 
$$
\t{B}_{ij}(e)=(0, \ldots, e_j \times e_i, 0, \ldots, e_i \times e_j , 0, \ldots )
$$
and hence $\t{B}_{ij}$  vanishes at $e$ and is tangent to $\t{M}_r$. The 
induced field on $\t{M}_r$ will 
be denoted $B_{ij}'$. Then $B_{ij}'$ projects to $B_{ij}$ on $M_r$. 
We note $dim T_e (\t N_r)=2n$, 
$dim T_e (\t M_r)=2n-2$ and $dim T_e (M_r)=2n-4$.

\begin{rem}
Since $e$ is a singular point of $M_r$ we have
$$
dim T_e(M_r) = 2n-4 > dim M_r = 2n-6.
$$
\end{rem}

We 
will first compute the linearization of $\t B_{ij}$ at $e$ in $\t{N}_r$ 
($e$ is a smooth point on $\t{N}_r$ so we use the first procedure) 
to obtain $\t A_{ij}\in End(T_e(\t{N}_r))$. 
Then $\t A_{ij}$ will preserve the subspace $T_e(\t{M}_r)\subset T_e(\t{N}_r)$ 
whence 
we obtain an induced element $A_{ij}'\in End(T_e(\t{M}_r))$. 
But there is an exact sequence
$$
V_e \to T_e(\t{M}_r) \to T_e(M_r)
$$
where $V_e= \{\del: \exists v\in \R^3 \hbox{~~such that~~} 
\del_i= e_i\times v, 1\le i\le n\}$, we note that $dim V_e=2$. 
We  will verify that 
$A_{ij}'(V_e)\subset V_e$ (in fact $A_{ij}(V_e)=0$). Hence $A_{ij}'$ 
will descend to $T_e(M_r)$. The resulting element of $End(T_e (M_r))$ will be $A_{ij}$, the 
linearization of $B_{ij}$ at $e$. 

Accordingly we begin by computing the linearization $\t A_{ij}$ of $\t B_{ij}$ on 
$T_e(\t{N}_r)$. Thus $\t A_{ij}$ will be $2n\times 2n$ matrix (instead of a 
$2n-4 \times 2n-4$ matrix). 
% and will be the linearization of $B_{ij}$ at $e\in M_r$. 

Another advantage in passing to $\t{N}_r$ is that $T_e(\t{N}_r)$ 
is now a direct sum of the tangent bundles of the factors
$$
T_e(\t{N}_r)= \oplus_{i=1}^n T_{e_i}(S^2(r_i)).
$$
The Riemannian connection on $\t{N}_r$ is a direct sum of the 
Riemannian connections on the summands. Thus we may write 
(for $\del\in T_e(\t{N}_r)$) 
$$
\t A_{ij}(\del)= (0,..., \nabla_{\del}(e_j \times e_i), 0,..., 
\nabla_{\del}(e_i\times e_j), 0,...).
$$
We will suppress the zeroes in the above row vectors henceforth. 

\begin{lem}
$$\t{A}_{ij}(\del)= (u\times \del_i, u\times \del_j) 
\left[ \begin{array}{cc} 
\eps_j r_j &  -\eps_j r_j\\
-\eps_i r_i &  \eps_i r_i 
\end{array}\right]$$
\end{lem}
\proof In the above formula for $\t{A}_{ij}(\del)$ we use the Riemannian 
connection $\nabla$ on $S^2$. We will compute using the flat connection 
$\bar\nabla$ on $T(\R^3)|S^2$ and then project back into $T(S^2)$ to get $\nabla$. 
We have
$$
\bar\nabla_{\del}(e_j \times e_i) = \del_j \times e_i + e_j \times \del_i
$$
$$
\bar\nabla_{\del}(e_i \times e_j) = \del_i \times e_j + e_i \times \del_j.
$$
Evaluating at $e$ we obtain
$$
\bar\nabla_{\del}(e_j \times e_i)|_e = \eps_i r_i \del_j\times u + 
\eps_j r_j u \times \del_i= \eps_j r_j u \times \del_i- \eps_i r_i u\times \del_j.
$$
Since the right-hand side is in $T_e(S^2)$ we have also
$$
\nabla_{\del}(e_j \times e_i)|_e = 
\eps_j r_j u \times \del_i- \eps_i r_i u\times \del_j.
$$
Finally $\nabla_{\del}(e_i \times e_j)|_e= -\nabla_{\del}(e_j \times e_i)|_e$ 
and the lemma follows. \qed

We now relate the action of $\p_n$  on $T_e(\t{N}_r)$ we have just computed 
to the action on $T_e(M_r)$. 
We recall that $\t{M}_r= \{e\in \t N_r : \sum_{i=1}^n e_i =0\}$ whence 
 $T_e(\t{M}_r)= \{\del\in T_e(\t N_r) : \sum_{i=1}^n \del_i =0\}$. 
  We have the 2-dimensional subspace $V_e$ of tangents to the $SO(3)$-orbit through 
$e$ described above. 
Thus we have a filtration $F_{\bullet}$ given by
$$
V_e\subset T_e(\t{M}_r) \subset T_e(\t{N}_r)
$$
and a canonical isomorphism 
$$
T_e(\t{M}_r)/V_e \cong T_e({M}_r). 
$$
We now show that $\p_n$ preserves the above filtration. 

\begin{lem}
\begin{enumerate}
\item $\p_n T_e(\t{N}_r)\subset T_e(\t{M}_r)$. 
\item $\p_n V_e =0$. 
\end{enumerate}
\end{lem}
\proof (1) is immediate. We prove (2). Suppose $\del\in V_e$. We claim 
$$
\del_j \times e_i + e_j \times \del_i=0, \quad 1\le i< j\le n . 
$$ 
Indeed, 
$$
\del_j \times e_i + e_j \times \del_i= (e_j \times v)\times e_i + 
e_j\times (e_i \times v)= (e_j \times v)\times e_i + (e_j \times e_i)
\times v + e_i\times(e_j\times v).$$
But $e$ is degenerate, so $e_i\times e_j=0$. \qed 
 
We collect our results in 

\begin{thm}
\begin{enumerate}
\item There is a $\p_n$-stable filtration 
$$
V_e\subset T_e(\t{M}_r)\subset T_e(\t{N}_r).
$$
\item $T_e(M_r)\cong  T_e(\t{M}_r)/V_e$.
\item There is an isomorphism
$$
\phi: T_e(\t{N}_r) \to T_u(S^2)\otimes \R^n
$$
such that $\phi \rho \phi^{-1} (X_{ij}) = 
ad u \otimes J_{ij}(\eps_i r_i, \eps_j r_j)$.
\item $\phi(T_e(\t{M}_r))= T_u(S^2)\otimes \R^n_0$ and 
$\phi(V_e)= T_u(S^2)\otimes \R v(\eps,r)$. 
Here $\R^n_0= \{(x_1,...,x_n): \sum_{i=1}^n x_i =0\}$ and 
$v(\eps,r)= (\eps_1 r_1,..., \eps_n r_n)$. 
\end{enumerate}
\end{thm}

Here $\R^n$ is realized as the space of row vectors with $n$ components.

\section{The action on the holomorphic tangent space.} 

The point of this section is that $T_e(\t{N}_r)$ has a $\p_n$-invariant 
almost complex structure that descends to  $T_e(M_r)$. We will compute 
the corresponding action of $\p_n$ on the holomorphic tangent space. 

Define an almost complex structure $J\in End (T_e(\t{N}_r))$ by 
$$
J(\del)= \eta \hbox{~~such that~~} \eta_i= u\times \del_i, 1\le i\le n\quad . 
$$
%The following lemma is immediate

\begin{lem}
\begin{enumerate}
\item $J$ is $\p_n$-invariant. 
\item The filtration $F_{\bullet}$ is invariant under $J$. 
\end{enumerate}
\end{lem}
\proof The first assertion is immediate. It is also clear that $T_t(\t M_r)$ is invariant under 
$J$. It remains to check that $V_e$ is invariant under $J$. Suppose $\del\in V_e$. 
Hence there exists $v\in \R^3$ such that $\del_i=\eps_i r_i u\times v$, $1\le i\le n$. Then 
$J\del_i= u\times (\eps_i r_i u \times v)=\eps_i r_i u\times (u\times v)$. Hence if we put 
$w=u\times v$ then 
$$
J \del_i = \eps_i r_r u\times w, 1\le i\le n \quad . 
$$
Therefore $J\del \in V_e$. $\qed$

\begin{rem}
The almost complex structure $J$ is {\bf not} the one induced by the 
complex structure on $\t{N}_r= \prod_{i=1}^n S^2(r_i)$. We have 
changed the complex structure on $S^2(r_i)$ to its conjugate for 
each $i$ such that $e_i$ is a back-track (i.e. $\eps_i=-1$).  
\end{rem}

We can decompose $T_e(\t{N}_r)\otimes \C$ into the $+i$-eigenspace of $J$ 
denoted by $T_e^{\eps}(\t{N}_r)$ and the $-i$-eigenspace denoted by
$T_e^{-\eps}(\t{N}_r)$.  Accordingly we have
$$
T_e^{\eps}(\t{N}_r)=\{ \del \in T_e(\t{N}_r)\otimes \C : u\times \del_j = \sqrt{-1}\del_j\}
$$
Similarly we denote the $+i$-eigenspaces of $J$ acting on 
$T_e(\t{M}_r)\otimes \C$ and $V_{e}\otimes \C$ by
$T_e^{\eps}(\t{M}_r)$ and $V^{\eps}_e$ respectively. We denote the 
quotient 
$T_e^{\eps}(\t{M}_r)/V^{\eps}_e$ by $T_e^{\eps}(M_r)$. Clearly the 
latter space is the $+i$-eigenspace of $J$ acting on $T_e({M}_r)\otimes \C$. 

Now we recall that we have an isomorphism
$$
\phi: T_e(\t N_r) \to T_u(S^2)\otimes \R^n
$$
complexifying we obtain
$$
\phi: T_e(\t N_r)\otimes \C \to T_u(S^2)\otimes_{\R} \C^n.
$$
We see that $\phi$ conjugates $J$ to $ad u \otimes 1$ and we have 
an induced isomorphism (again denoted by $\phi$)
$$
\phi: T_e^{\eps}(\t N_r) \to T_u^{1,0}(S^2)\otimes_\C \C^n.
$$
Under $\phi$ the action of $X_{ij}$ transforms to 
$\sqrt{-1} I\otimes J_{ij}(\eps_i r_i, \eps_j r_j)$. We note that 
\sloppy{$dim_\C T_u^{1,0}
(S^2)=1$} and we obtain a canonical isomorphism
$$
\psi: T_e^{\eps}(\t N_r) \to\C^n.
$$
This isomorphism has the property:
$$
\psi(T_e^{\eps}(\t M_r))=\C^n_0,\quad \psi(V^{\eps}_e)= \C v(\eps, r).
$$
We have completed our computation of the action of $\p_n$. 

\begin{thm}
\begin{enumerate}
\item There is a canonical isomorphism $\psi: T_e^{\eps}(\t N_r) \to\C^n$
\item $\psi$ induces the action of $X_{ij}\in \p_n$ on $\C^n$ by 
$\sqrt{-1} J_{ij}(\eps_i r_i, \eps_j r_j)$. 
\item $\C^n$ admits a $\p_n$-invariant filtration by 
$\psi(T_e^{\eps}(\t M_r))=\C^n_0$, $\psi(V^{\eps}_e)= \C v(\eps, r)$. 
\item There is an $\p_n$-invariant complex structure $J$ on $T_e(M_r)$. 
The induced action of $\p_n$ on the $+i$-eigenspace of $J$ in 
$T_e(M_r)\otimes \C$ corresponds to the action of $\p_n$ on 
the quotient $\C^n_0/\C v(\eps, r)$. 
\end{enumerate}
\end{thm}

Here $\C^n$ is realized as the space of row vectors with $n$ components.

\section{The associated hypergeometric equation.} 

As discussed in the introduction we use the linear operators 
$A_{ij}\in End(\C^n)$ to obtain a flat holomorphic connection 
$\nabla$ on the trivial $T^{\eps}_e(\t{N}_r)$-bundle $\e$  
over $\M= \C^n_*$. The connection form $\om$ of $\nabla$ is
$$
\om = \sum_{1\le i< j\le n} \frac{dz_i - dz_j}{z_i -z_j} \otimes A_{ij}.
$$
A (multivalued) holomorphic section of $\e$ corresponds to a {\em row} vector  
$F=(F_1,...,F_n)$ of (multivalued) holomorphic functions. The 
hypergeometric equation comes from the condition that $F$ be parallel 
 for the connection $\nabla$:
$$
dF =  F\om 
$$
or equivalently
\begin{equation}
\label{equ}%[(*)]
d F_i = \sum_{j, j\ne i} (\la_j F_i -\la_i F_j)\frac{dz_i - dz_j}{z_i -z_j} 
\end{equation}
with $\la_j= \sqrt{-1}\eps_j r_j$. We will refer to (\ref{equ}) as the 
{\em hypergeometric equation}. 

We observe that the operators $A_{ij}$ leave invariant the subspace $\C^n_0$ and annihilate the line 
$V_{\la}= \C (\la_1,..., \la_n)$. We obtain a diagram of flat bundles over $\C^n_*$:
$$
\begin{array}{ccc}
\C^n_* \times \C^n_0 & \lra & \C^n_* \times \C^n\\
\downarrow & ~ & ~ \\
\C^n_* \times \C^n_0 /V_{\la} & ~ & ~
\end{array}
$$ 

The monodromies of these bundles will be the representations of $P_n$ 
corresponding to the actions of ${\cal P}_n$ on 
$T_e(\t{N}_r)$, $T_e(\t{M}_r)$, $T_e({M}_r)$.

\section{Solving the hypergeometric equation by hypergeometric integrals.} 

Let $\la_1,..., \la_n$ be complex numbers with $\la_j\notin \Z, 1\le j\le n$. 
Let $(\xi, z_1,..., z_n)\in (\C^{n+1})_*$ and $\Phi(\xi, z_1,..., z_n)$ 
be the hypergeometric integrand
$$
\Phi(\xi, z_1,..., z_n):= (\xi- z_1)^{\la_1} \ldots (\xi - z_n)^{\la_n} .
$$ 
Let $\chi:=\chi_\la : \F_n \to \C^*$ be the character defined by 
$\chi(\ga_j)=\exp(2\pi \sqrt{-1}\la_j)$, $1\le j\le n$. Recall that   
$\{\ga_1,...,\ga_n\}$ is a generating set for $\F_n$, the free group of rank $n$. 
Here we identify $\F_n$ with the fundamental group $\pi_1(M, b)$, where 
$M= \C- \{z_1,..., z_n\}$, so that the conjugacy class of $\ga_j$ 
is represented by a sufficiently small loop 
which goes once around $z_j$ in the counterclockwise direction. 
Note that $\chi(\ga_j)\ne 1$, $1\le j\le n$.  
For any character $\chi: \F_n\to \C^*$ we let $L_\chi$ 
be the local system over $M$ given by 
$$
L_\chi = \t M \times \C /((x,z)\sim (\ga x, \chi(\ga)z)) .
$$  
We define a multivalued parallel section $\si$ of $L_\chi$ by $\si(x)= [x,1]$ 
(where $[x,z]$ denotes the equivalence class of $(x,z)$). Note that the lift 
of $\si$ to the universal cover satisfies
$$
\si(\ga x)= [\ga x, 1] = [x, \chi(\ga)^{-1}] = \chi(\ga)^{-1} \si(x).
$$
The following lemma is obvious:

\begin{lem}
The $L_\chi$-valued 1-forms $\zeta_j$, $1\le j\le n$, defined by
$$
\zeta_j(\xi)= (\xi- z_1)^{\la_1}\ldots (\xi- z_n)^{\la_n}\frac{d\xi}{\xi- z_j} 
\otimes\si
$$
are single-valued on $M$. 
\end{lem}
Hence $\zeta_j$ gives rise to a class $[\zeta_j]$ in the de Rham cohomology group $\H^1_{dR}(M, L_\chi)$. 

Let $\ga\in \H_1(M, L_{\chi^{-1}})$. Let $G_j$ be the Kronecker pairing $\< \zeta_j, \ga\>$ considered as a 
function of $z_1,..., z_n$. This Kronecker pairing is traditionally represented as an integral. 
To make this precise let $\ga=\sum_{i=1}^k a_i \otimes \tau_i$, where each $a_i$, $1\le i\le k$, is a  1-simplex 
and $\tau_i$ is a parallel section of $\L^{-1}|a_i$. Then $\<\zeta_j, \ga\>$ is given by 
$$
G_j(z_1,..,z_n)= \sum_{i=1}^k \int_{a_i} (\xi- z_1)^{\la_1}\ldots (\xi- z_n)^{\la_n}\<\si, \tau_i\> 
\frac{d\xi}{\xi- z_j} \quad . 
$$
We will use the following more economical notation:    
$$
G_j(z_1,..,z_n)= \int_{\ga} (\xi- z_1)^{\la_1}\ldots (\xi- z_n)^{\la_n}\frac{d\xi}{\xi- z_j}\otimes\si \quad . 
$$

Now we let $z=(z_1,...,z_n)$ vary. Let $\pi: \C^{n+1}_* \to \C^{n}_*$ be the map that forgets the first component. 
Then $\pi^{-1}(z)$ is isomorphic to $\C - \{z_1,...,z_n\}$. By \cite[3.13]{DM}, the flat line bundle 
$L_\chi$ on $\pi^{-1}(z)$ is the restriction of a flat line bundle $\t L_\chi$ on $\C^{n+1}_*$. As $z$ varies, 
the forms $\zeta_1,...,\zeta_n$ give rise to {\em relative holomorphic 1-forms} on $\C^{n+1}_*$ with  
coefficients in $\t L_\chi$. We recall that a relative holomorphic form on the total space $E$ of a 
holomorphic fiber bundle $p: E\to B$ is an element of the quotient differential 
graded algebra
$$
\Om^\bul (E)/ (p^* \Om^\bul (B)^+). 
$$
Here $\Om^q$ denotes the holomorphic $q$-forms and $(p^* \Om^\bul (B)^+)$ denotes the differential ideal 
in $\Om^\bul (E)$ generated by the pull-backs to $E$ of holomorphic forms on $B$ of positive degree. A 
relative holomorphic $q$-form $\eta$ is relatively closed if $d\eta$ is in the above ideal. The forms 
$\zeta_1,...,\zeta_n$ are relatively closed, hence they induce holomorphic sections 
$[\zeta_1],...,[\zeta_n]$ of the vector bundle ${\mathcal H}^1$ over $\C^{n}_*$ with fiber over $z$ given by 
$$
\H^1(\pi^{-1}(z), \t L_\chi |\pi^{-1}(z)).
$$
Precisely, $[\zeta_i](z)$ is the class of the 1-form $\zeta_i(z)$ on $\pi^{-1}(z)$ in the above cohomology group. 
The bundle ${\mathcal H}^1$ has a flat connection, the {\em Gauss-Manin} connection, whose definition we now 
recall. Note first that a local trivialization of $\pi$ induces a local 
trivialization of ${\mathcal H}^1$. Then a smooth section of ${\mathcal H}^1$ is parallel for the 
Gauss-Manin connection if it is constant when expressed in terms of all such induced local trivializations. The bundle 
${\mathcal H}_1$ of the first homology groups with coefficients in $\t L_{\chi^{-1}}$ admits an analogous flat connection. 
Now let $p: \t\C^n_* \to \C^n_*$ denote the universal cover of $\C^n_*$. We obtain a pull-back fiber bundle 
$\t\pi: E\to \t\C^n_*$ of $n$-punctured complex lines over $\t\C^n_*$ and pull-back flat vector bundles 
$\t{\mathcal H}^1$ and $\t{\mathcal H}_1$. Choose a base-point $z^0=(z_1^0,...,z_n^0)$ in $\C^n_*$. We use $M$ to denote 
$\C - \{z_1^0,...,z_n^0\}$ henceforth. Choose a base-point $\t z^0$ in  $\t\C^n_*$ lying over $z^0$. 
We may identify the fiber of $\t{\mathcal H}_1$ over $\t z^0$ with $\H_1(M, L_{\chi^{-1}})$. Hence 
given $\ga\in \H^1(M, L_{\chi^{-1}})$ there is a unique parallel section $\t\ga$ of 
$\t{\mathcal H}_1$ such that $\t\ga(\t z^0)=\ga$. We can now define a global holomorphic 
function $G_j(z)$ on $\t\C^n_*$ by
$$
G_j(z)= \int_{\t\ga} (\xi -z_1)^{\la_1} \ldots  (\xi -z_n)^{\la_n} \frac{d\xi}{\xi -z_j} \otimes \si. 
$$
Here we have used the same notation for corresponding (under pull-back) objects on $\C^n_*$ and 
$\t\C^n_*$. We may also write 
$$
G_j(z)= \< [\zeta_j(z)], \t\ga\>
$$
where $\< , \>$ is the fiberwise pairing between $\t{\mathcal H}^1$ and $\t{\mathcal H}_1$. We have

\begin{lem}
\label{gm}
$$
dG_i(z)= \sum_{j=1}^n (\int_\ga \frac{\D}{\D z_j} (\frac{\Phi}{\xi- z_i}) d\xi \otimes \si )dz_j. 
$$
\end{lem}
\proof We have
$$
dG_i(z)= \< \nabla [\zeta_i(z)], \t\ga\>
$$
where $\nabla$ is the Gauss-Manin connection. We will need another formula for the 
Gauss-Manin connection, see \cite{KO} or Remark \ref{gmr} below. 
Before stating the formula we need more notation. 
Let $F^q \Om^q(E)$ denote the subspace of holomorphic $q$-forms on $E$ that 
are multiples of pull-backs of  $q$-forms from the base $\t\C^n_*$ by elements 
of ${\mathcal O}(E)$.  Then we have a canonical isomorphism (because the fibers of $\t\pi$ 
have complex dimension $1$)
$$
\frac{\Om^2(E)}{d F^1\Om^1(E) + F^2\Om^2(E)} 
\cong \Om^1(\t\C^n_*, \t{\mathcal H}^1).
$$
Now the formula for $\nabla$ is
$$
\nabla [\zeta_i] = [d\zeta_i]. 
$$
Here $d\zeta_i$ denotes the exterior differential of $\zeta_i$ where $\zeta_i$ is considered as a 1-form on $E$ (modulo  
$F^1\Om^1(E) $) with values in the line bundle $p^* \t L_\chi$. The symbol 
$[d\zeta_i]$ denotes the class of $d\zeta_i$ modulo $d F^1\Om^1(E) + F^2\Om^2(E)$. The lemma follows from 
the formula
$$
d\zeta_i\equiv \sum_{j=1}^n  \frac{\D}{\D z_j} (\frac{\Phi}{\xi- z_i}) dz_j \wedge d\xi \otimes \si 
$$
together with the observation that integration over $\t\ga$ factors through $[ \ ]$. \qed

\begin{rem}\label{gmr}
The above formula for $\nabla$ can be proved as follows. First note that the formula does indeed 
define a connection, to be denoted $\nabla'$ on ${\mathcal H}^1$. To show that $\nabla$ and $\nabla'$ 
agree it suffices to show they agree locally. Since they are both invariantly defined it suffices 
to prove that they agree on trivial bundles. But it is clear that in this case 
a section of ${\mathcal H}^1$ is parallel for $\nabla'$ if and only if it is constant. 
\end{rem}

The proof of the next lemma is a modification of 
\cite[Proposition 2.2.2]{K1}. 

\begin{lem}
The functions $G=(G_1,..., G_n)$ satisfy
$$
dG_i = \sum_{j, j\ne i} (\la_j G_i - \la_j G_j) \frac{dz_i - dz_j}{z_i-z_j} 
\otimes \si  \quad or \quad dG^T= \om G^T \quad . 
$$
\end{lem}
\proof We will drop the $\otimes \si$ for the course of the proof: 
$$
G_i(z)= \int_\ga \Phi \frac{d\xi}{\xi- z_i}. 
$$
Whence by Lemma \ref{gm} 
$$
dG_i =  -\sum_{j=1}^n [\int_{\ga} \la_j\Phi (\xi-z_j)^{-1} (\xi - z_i)^{-1} 
d\xi] dz_j   + [\int_{\ga}\Phi(\xi- z_i)^{-2} d\xi]dz_i$$
$$ 
= -\sum_{j\ne i} [\int_{\ga} \la_j\Phi (\xi-z_j)^{-1} (\xi - z_i)^{-1} d\xi] dz_j
-  [\int_{\ga}(\la_i-1)\Phi(\xi- z_i)^{-2} d\xi]dz_i .
$$
We simplify the first term using
$$
\frac{1}{\xi- z_i} \cdot \frac{1}{\xi - z_j} = \frac{1}{z_i -z_j} 
(\frac{1}{\xi- z_i} - \frac{1}{\xi - z_j})
$$
to obtain 
$$
= -\sum_{j\ne i} \frac{\la_j}{z_i -z_j} 
[ \int_{\ga}\Phi \frac{d\xi}{\xi-z_i}   
- \int_{\ga}\Phi \frac{d\xi}{\xi-z_j}]dz_j 
- [\int_{\ga}(\la_i-1)\Phi(\xi- z_i)^{-2} d\xi]dz_i = $$
$$
= -\sum_{j\ne i} \frac{\la_j G_i}{z_i-z_j} dz_j + 
\sum_{j\ne i} \frac{\la_j G_j}{z_i-z_j} dz_j - [\int_{\ga}(\la_i-1)\Phi(\xi- z_i)^{-2} d\xi]dz_i .$$
Now we have
$$
d(\Phi(\xi- z_i)^{-1}) = (\la_i -1) \Phi(\xi - z_i)^{-2} d\xi + 
\sum_{j\ne i}\la_j \Phi(\xi- z_i)^{-1}(\xi- z_j)^{-1} d\xi .$$
Thus by Stokes' Theorem
$$
-\int_{\ga}(\la_i -1)\Phi(\xi-z_i)^{-2}d\xi = \int_{\ga} \sum_{j\ne i}\la_j \Phi(\xi- z_i)^{-1}(\xi- z_j)^{-1} d\xi =$$
$$
\int_{\ga} \sum_{j\ne i}\la_j\Phi \frac{1}{z_i -z_j}
(\frac{1}{\xi- z_i} - \frac{1}{\xi- z_j}) d\xi = 
$$
$$
\sum_{j\ne i}\frac{\la_j}{z_i-z_j} G_i - \sum_{j\ne i}\frac{\la_j}{z_i-z_j} G_j
$$
hence 
$$
-[\int_{\ga}(\la_i -1)\Phi(\xi-z_i)^{-2}d\xi]dz_i= 
\sum_{j\ne i}\frac{dz_i}{z_i-z_j} (\la_j G_i - \la_j G_j)\quad .
$$
We obtain 
$$
dG_i= \sum_{j\ne i}\frac{dz_i-dz_j}{z_i-z_j} (\la_j G_i - \la_j G_j)\ . \qed
$$

\begin{rem}
The simplification using Stokes' Theorem above is equivalent to observing that
$$
\Phi (\xi - z_i)^{-1} dz_i \otimes \si \in F^1 \Om^1(E), 1\le i\le n, 
$$
and we work modulo $dF^1 \Om^1(E)$ in computing $\nabla$. 
\end{rem}

We now define $F_i:= \la_i G_i$, $1\le i\le n$. 

\begin{lem}
$F= (F_1,..., F_n)$ is a solution of the hypergeometric equation (\ref{equ}). 
\end{lem}
\proof 
$$
dF_i = \la_i dG_i = \sum_{j\ne i}\frac{dz_i-dz_j}{z_i-z_j} 
(\la_i\la_j G_i - \la_i\la_j G_j)=$$
$$ 
=\sum_{j\ne i} \la_j(\la_i G_i) - \la_i (\la_jG_j) 
\frac{dz_i-dz_j}{z_i-z_j}= 
$$
$$
=\sum_{j\ne i} (\la_j F_i -\la_i F_j) \frac{dz_i-dz_j}{z_i-z_j} \ . \qed
$$

We have proved

\begin{thm}
Let $\ga$ be an element of $\H_1(M, L_{\chi^{-1}})$ 
and $\si$ a flat multivalued section of $L_{\chi}$. 
For $\la=(\la_1,..., \la_n)\in \C^n$ define a holomorphic function on $\t\C^n_*$ by
$$
F_i:= \la_i \int_{\t\ga} (\xi-z_1)^{\la_1}\ldots (\xi- z_n)^{\la_n}\frac{d\xi}{\xi -z_i}\otimes\si .
$$
Then $F=(F_1,...,F_n)$ is a solution of the hypergeometric equation.
\end{thm}

\section{The monodromy representation of the hypergeometric equation and the action on homology.}

We have seen that for $\ga\in \H_1(M, L_{\chi^{-1}})$ 
we obtain a solution $S= (F_1,..., F_n)$ of the hypergeometric equation by the formula
$$
F_i:= \la_i \int_{\t\ga} (\xi-z_1)^{\la_1}\ldots (\xi- z_n)^{\la_n}\frac{d\xi}{\xi -z_i}\otimes \si .
$$
It is important to recall that $\sum_{j=1}^n\la_j=0$. The differential forms
$$
\eta_j= \la_j (\xi-z_1)^{\la_1}\ldots (\xi- z_n)^{\la_n}\frac{d\xi}{\xi -z_i}\otimes \si 
$$
are de Rham representatives of the cohomology classes $[\eta_j], 1\le j\le n,$ in 
$\H^1(M, L_{\chi^{-1}})$. Note that 
$$
d((\xi -z_1)^{\la_1}\ldots (\xi- z_n)^{\la_n}\otimes \si)= \eta_1 +\ldots \eta_n
$$
hence we have the relation
\begin{equation}
\label{eq1}
 [\eta_1]+\ldots +[\eta_n]=0
\end{equation}

\begin{lem}
The span of the cohomology classes 
$[\eta_j], 1\le j\le n$, has dimension $n-1$.  
\end{lem}
\proof First since $\sum_{j=1}^n \la_j=0$ we have $\chi(\ga_1 \ga_2 \ldots 
\ga_n)=1$. Thus  
$L_{\chi}$ extends to a flat line bundle over $\C\P^1 - \{z_1,\ldots , z_n\}$. 
Also, $\eta_j$ extends meromorphically over infinity with a simple pole at 
infinity. 

Next we extend the flat line bundle $L_\chi$ to a holomorphic line bundle 
$\L^{hol}$ on $\C\P^1$ so that $(\xi-z_j)^{\la_j}\otimes\si$ is a local basis 
around $z_j$. Then  $(\xi-z_1)^{\la_1}\ldots (\xi- z_n)^{\la_n}\otimes \si$ 
is a holomorphic section of $\L^{hol}$ which has no zeroes or poles. 

We can now prove the lemma. We have a flat line bundle $L_\chi$ over 
$M$ (with trivial monodromy around $\infty$). 
The argument of \cite[\S 2.7]{DM} proves that we can compute 
the group $\H^1(M, L_\chi)$ as the 1-st cohomology group of the complex 
$(\Om^{\bullet}(\C\P^1, *D, L_\chi), d)$ of holomorphic 
$L_\chi$-valued forms on $M$ which have at worst poles at 
$z_1,..., z_n, \infty$. Here the (additive) divisor $D$ 
is defined by $D=z_1 + ... +z_n + \infty$. Now  
$\eta_j \in \Om^{1}(\C\P^1, *D, L_\chi)$ and 
$$
\Om^{0}(\C\P^1, *D, L_\chi)= \{ f\Phi\otimes \si : \hbox{~~so that~~} f 
\hbox{~~has at worst poles at~~} D \}.
$$
First note that $Span(\eta_1,..., \eta_n)\subset \Om^{1}(\C\P^1, *D, L_\chi)$ 
has dimension $n$ since the forms $\eta_j$ have singularities at distinct points of $\C$. 

Suppose that there exists $f\Phi\otimes \si\in \Om^{0}(\C\P^1, *D, L_\chi)$ and $c_1,..., c_n$ such that
$$
d(f \Phi \otimes \si )= c_1 \eta_1 +...+ c_n \eta_n \quad .
$$
We claim that $f$ cannot have any poles. Indeed, assume $f$ has a pole of order $k\ge 1$ at $z_i$. Then 
$$
f(\xi)= \frac{c}{(\xi - z_i)^{k}} + \ldots 
$$
We are assuming
$$
df \Phi + f d\Phi = \sum_{i=1}^n c_i \eta_i
$$
or
\begin{equation}
\label{eq2}
df \Phi + (f \sum_{i=1}^n \frac{\la_i}{\xi - z_i} d\xi ) \Phi = \sum_{i=1}^n c_i \eta_i .
\end{equation}
Equating the coefficients of $(\xi - z_i)^{-k-1}$ in the equation (\ref{eq2}) from each side we obtain
$-k c + \la_i c=0$, or $\la_i=k$. This contradicts the assumption that each $\la_i$ is pure imaginary. 
It remains to check that $f$ is not a polynomial. Assume $f$ has a pole of order $k\ge 1$ at $\infty$, whence 
$f(\xi)= a_0 + a_1 \xi + ... + a_k \xi^k$. We equate the coefficients at $\xi^{k-1}d\xi$ on each side 
of (\ref{eq2}) to obtain
$k a_k + (\sum_{i=1}^n \la_i) a_k =0$ or $ka_k=0$. This contradiction proves the claim. 
 Hence $f\equiv c$ and hence 
$$
df= c\sum_{i=1}^n \eta_i
$$
which means that the dimension of the subspace of coboundaries in 
$Span(\eta_1,..., \eta_n)$ is $1$.  \qed

In the group cohomology computations that follow $\ga_1,..., \ga_n$ will be a 
generating set of $\F_n$ and $b_1,..., b_n$ will be its image under abelianization in 
$\Z^n$. Here the loop representing $\ga_i$ is obtained by connecting the small circle $a_i$ 
going around $z_i$ to the base-point $*\in \C - \{ z_1,..., z_n\}$. 
We recall that $P_n$ acts on $\F_n$ preserving the conjugacy classes of the 
generators $\ga_j$. Hence the induced action on $\Z^n$ is trivial and $P_n$ fixes any 
character $\chi: \F_n\to \C^*$. Hence $P_n$ acts on $\H^1(\F_n,\C_{\chi})$. 
Here we let $\C_{\chi}$ denote the 1-dimensional space on which $\F_n$ acts via $\chi$. 
We next need

\begin{lem}
Suppose that $\chi: \F_n\to \C^*$ satisfies $\chi(\ga_i)\ne 1$ for all $i$. 
Then 
\sloppy{$dim_{\C} \H^1(\F_n,\C_{\chi})=n-1$.} 
\end{lem}
\proof The Euler characteristic $E(\F_n, \C_1)=1-n$. Hence $E(\F_n, \C_\chi)=1-n$. 
On the other hand, $\H^0(\F_n, \C_{\chi})=0$. \qed

\begin{cor}
$dim_{\C} \H_1(M, L_{\chi^{-1}})=n-1$ and the classes $[\eta_1],..., [\eta_{n-1}]$ form a basis for 
 $\H^1(M, L_{\chi})$ 
\end{cor}

We can construct an explicit basis $w_1,...,w_{n-1}$ for 
$\H_1(M, \L^{-1})$ following \cite[\S 2]{DM} as follows. We write 
$w_i= \ga_i \otimes \si_i + \ga_{i+1}\otimes \si_{i+1}$, where $\si_i, \si_{i+1}$ are  
multivalued flat sections along $\ga_i, \ga_{i+1}$ respectively and the jump experienced by 
$\si_i$ (at the base-point) after parallel translating along $\ga_i$ cancels 
that of $\si_{i+1}$ along $\ga_{i+1}$.  

Define flat sections $S_i$, $1\le i\le n-1$, of $\t\C^n_*\times \C^n_0$ by 
$$
S_i:= (S_{i1},..., S_{in}), \hbox{~~where~~} S_{ij}= \la_j \int_{\t w_i}\eta_j.
$$
We see then that $S_1,..., S_{n-1}$ are  multivalued parallel sections of 
$\C^n_*\times \C^n_0$. %These section span the space of parallel sections 
%of $\C^n_*\times \C^n_0$. To get a basis of the 
%space of parallel sections of $\C^n_*\times \C^n$ we add to this collection the 
%constant section $S_n=(\la_1,...,\la_n)$. 

The desired representation $\rho: P_n\to Aut(\C^n_0)$ 
is obtained by parallel translation of $S_1,..., S_{n-1}$ along loops in $\C^n_*$. 
The resulting automorphisms leave invariant the 
line $\C\la$ where $\la=(\la_1,...,\la_n)$.

Before stating the main result of this section we need to define a special class $w_{\infty}$ in 
$\H_1(M, L_{\chi}^{-1})$. 
Let $a_{\infty}\subset \C$ be a circle whose interior contains all the punctures $z_1,...,z_n$. Since 
$\la_1 + \ldots \la_n=0$, the monodromy of $L_{\chi}^{-1}$ around $a_{\infty}$ is trivial. 
Hence there is a nonzero parallel section $\si^{\vee}$ of $L_{\chi}^{-1}|a_{\infty}$. We let 
$w_{\infty}$ be the homology class represented by $a_{\infty}\otimes \si^{\vee}$. 

Let $\tau: P_n \to Aut \H_1 (M, L_{\chi}^{-1})$ be the homomorphism induced by the inclusion 
$P_n \subset Aut(\F_n)$ (recall that $P_n$ acts trivially on the sheaf of parallel sections of $L_{\chi}^{-1}$). 

\begin{lem}
\label{nitty}
(1) $\int_{w_{\infty}} \eta_i= -\la_i$, in particular $w_{\infty}\ne 0$. 

(2) The class $w_{\infty}$ is fixed by $P_n$. 
\end{lem}
\proof To prove (1) we apply the residue theorem and note that 
$$
\Phi(\xi, z)|_{\xi=\infty} = 1
$$
and the residue of $(\xi - z_i)^{-1}d\xi$ at $\xi=\infty$ is $-1$. To verify (2) 
we identify $P_n$ with a subgroup of the mapping class group of $M$. 
Then we choose representatives for the elements of $P_n$ so that they act by the identity on the 
closure of the exterior of the circle $a_{\infty}$. $\qed$ 

We now have

\begin{thm}
(i) The monodromy representation of the flat bundle $\C^n_* \times \C^n_0$ is equivalent to $\tau$. 

(ii) Under the above equivalence the invariant line $V_{\la}\subset \C^n_0$ corresponds to the line 
$\C w_\infty \subset \H_1(M, \L^{-1}_{\chi})$. 

(iii) We obtain an induced equivalence of the monodromy representation of $\C^n_* \times \C^n_0/V_\la$ and the 
induced action of $P_n$ on $\H_1(\C\P^1-\{z_1,...,z_n\}, L_{\chi}^{-1})$. 
\end{thm}
\proof We have an isomorphism $\Psi$ from $\H_1(M, L_{\chi}^{-1})$ onto the space of parallel 
sections on $\t\C^n_*\times \C^n_0$ given by $\Psi(w)=S_w$ where
$$
S_w=(\int_{\t w} \eta_1,..., \int_{\t w} \eta_n)= (\<[\eta_1], \t w\>, \ldots, \<[\eta_n], \t w\>). 
$$
We claim that $\Psi$ intertwines the representations $\tau$ and $\rho$ (see above) of $P_n$. The monodromy representation 
$\rho: P_n \to Aut(\C^n_0)$ is defined by
$$
S_w(g^{-1} z)= S_w(z) \rho(g). 
$$
In order to go further we will need to lift the $P_n$ action on $\t\C^n_*$ to the total space 
of $\t\pi: E \to \t\C^n_*$. We note that from the fiber bundle 
$\pi: \C^{n+1}_* \to \C^{n}_*$ we get an exact sequence $\F_n\to P_{n+1}\to P_n$. We may split 
this sequence by mapping $P_n$ to the subgroup of $P_{n+1}$ which consists of those elements 
that do not involve the first string of a braid -- recall that $\pi$ forgets the first point.  
Let   $\t\C^{n+1}_*$ be the universal cover of $\C^{n+1}_*$. Then $P_{n+1}$ acts on 
$\t\C^{n+1}_*$. But $E= \t\C^{n+1}_*/\F_n$, whence $P_n=P_{n+1}/\F_n$ acts on $E$ as the group of 
deck transformations of the cover $E\to \C^{n+1}_*$, and we obtain 
the required 
lift $\t g$ of elements $g\in P_n$ to $Aut(E)$. 
We now can give a formula for the monodromy representation 
$\tau$, namely
$$
\t w (gz)= \t g_* \tau(g)^{-1} \t w(z)
$$ 
or
$$
\t w (g^{-1}z)= \t g_*^{-1} \tau(g) \t w(z). 
$$
We can now prove the claim. Observe that since $\eta_i$ is an invariantly defined 1-form 
with values in $L_\chi$ on $\C^{n+1}_*$ we have
$$
\eta_i(g z)= (\t g^{-1})^* \eta_i (z)
$$
or
$$
\eta_i(g^{-1} z)= (\t g)^* \eta_i (z). 
$$
Hence 
$$
S_w(z) \rho(g)= S_w(g^{-1} z)= (\int_{\t w (g^{-1}z)}  
\eta_1 (g^{-1}z),..., \int_{\t w (g^{-1}z)}  \eta_n (g^{-1}z))=
$$
$$
(\int_{g^{-1}_* \tau(g)\t w(z)} \t g^* \eta_1 (z),..., 
\int_{g^{-1}_* \tau(g)\t w(z)} \t g^*  \eta_n (z))=
$$
$$
(\int_{\tau(g) \t w (z)}  \eta_1 (z),..., \int_{\tau(g) \t w (z)}  \eta_n (z))
$$
and the claim is proved. Hence (i) follows.

To verify (ii) it suffices to observe that $S_{w_\infty} = (-\la_1,\ldots, -\la_n)$, which follows from 
Lemma \ref{nitty}. From (i) and (ii) we deduce that the monodromy 
representation of $\nabla$ on $\C^n/V_\la$ is equivalent to the action of $P_n$ on 
$\H_1(M,\L_\chi^{-1})/\C w_{\infty}$. But it is clear from the  exact sequence 
of the pair $(M,\C\P^1 -\{z_1,...,z_n\})$ that we have a natural isomorphism 
$\H_1(M,\L_\chi^{-1})/\C w_{\infty}\cong \H_1(\C\P^1-\{z_1,...,z_n\},\L_\chi^{-1})$. $\qed$

\begin{rem}
Since we have seen that $T_e(\t{M}_r)$ contains an invariant line, the corresponding representation of $P_n$ 
must be on $\H_1(M, \L^{-1})$, not on $\H^1(M, \L)$ (the latter has an invariant hyperplane). 
\end{rem}

\section{The Gassner Representation.}
\label{Gassner}

We will follow \cite{Bi} and \cite{Mo} for our treatment of the Gassner representation. We begin with 
a quick review of the Fox  calculus. 

Let $G$ be a finitely generated group and $M$ a $G$-module. Let  $\C[G]$ be the 
group ring. 

\begin{defn}
A  {\bf derivation} $D: \C[G]\to  M$ is a $\C$-linear map satisfying
$$
D(fh)= (D(f))\eps(h) + f D(h)
$$
where $\eps: \C[G]\to \C$  is the augmentation. We let $Der(G,M)$ denote the space of derivations. 
\end{defn}

\begin{rem}
The restriction of each derivation $D$ to $G$ is a 1-cocycle $\del\in Z^1(G,M)$. Conversely, given a 
1-cocycle $\del\in Z^1(G,M)$ we define  a  derivation $D$ by
$$
D(\sum_{i=1}^n c_i g_i)= \sum_{i=1}^n c_i \del(g_i).
$$ 
Thus $Der(G,M)$ and $Z^1(G,M)$ are canonically isomorphic. We will identify them henceforth. 
\end{rem}

In the case $G$ is the free group $\F_n$ on the generators $\{x_1,.., x_n\}$ there is a unique 
derivation $\frac{\D}{\D x_i}\in Der(\F_n, \C[F_n])$ given by
$$
\frac{\D}{\D x_i}(x_j)= \del_{ij}, 1\le i, j\le n.
$$
Then $Der(\F_n, \C[F_n])$ is free over $\C[F_n]$ with the  basis 
 $\frac{\D}{\D x_1}, \ldots , \frac{\D}{\D x_n}$. Note that the projection 
$p: \F_n \to \H_1(\F_n)\cong \Z^n$ induces a ring-homomorphism 
$p:  \C[\F_n]\to \C[\H_1(\F_n)]$ and a push-forward map on derivations
$$
p_*: Der(\F_n, \C[\F_n]) \to Der(\F_n, \C[\H_1(\F_n)]).
$$
We may identify $\C[\H_1(\F_n)]$ with the $\C$-algebra 
$\L$ of Laurent polynomials in $t_1,..., t_n$. 
The space $Der(\F_n, \L)$ is free over $\L$ with the basis 
$p_*\frac{\D}{\D x_1}, \ldots , p_*\frac{\D}{\D x_n}$. 
We will drop $p_*$ henceforth. 

The main point in the construction of the Gassner representation is 
that there is a homomorphism 
$\si: P_n \hook Aut(\F_n)$. This homomorphism is described in 
terms of formulas in \cite[Corollary 1.8.3]{Bi}. There is an 
elementary description of $\si$ 
in terms of ``pushing a loop along the braid'', see 
\cite[Page 87]{Mo}. In both cases the action 
of $P_n$ on $\F_n$ is a {\em right} action, i.e. there is $\bar\si$ such that 
$\bar\si(p_1 p_2)= \bar\si(p_2)\bar\si(p_1)$. Therefore, the homomorphism $\si$ is 
actually given by $\si(p):= \bar\si(p^{-1})$. Next we note that 
we have an action of $P_n$ on  
$Der(\F_n, \L)$: 
$$
g\cdot D(x)= D(\si(g)^{-1} x).
$$
Since $P_n$ acts trivially on $\L$, $g\cdot D$ is still a derivation and the operator $g\cdot$ is $\L$-linear. 

\begin{rem}
In \cite{Bi} and \cite{Mo} the  action of $P_n$ on 
 $Der(\F_n, \L)$ is defined by $g* D(x)= D(\bar\si(g) x)$. But 
$\bar\si(g)= \si(g)^{-1}$ and hence $g\cdot D= g* D$. The composition of two right actions is a homomorphism!
\end{rem}

We can now define the Gassner representation. 

\begin{defn}
The Gassner representation $\rho: P_n \to Aut_\L  (  Der(\F_n, \L))$ assigns to each $g\in P_n$ the 
operator $g\cdot$ on $Der(\F_n, \L)$, where $Der(\F_n, \L)$ 
is considered as a free $\L$-module of rank $n$.  
\end{defn}

It is traditional to represent $\rho(g)$ as an element $(a_{ij})$ of $GL_n(\L)$ using the basis  
$\frac{\D}{\D x_1}, \ldots , \frac{\D}{\D x_n}$, see 
\cite[Page 119]{Bi}, \cite[Page 194]{Mo}:
$$
a_{ij}= \frac{\D}{\D x_j}\bar\si(g) x_i |_{x_i =t_i} .
$$

The Gassner representation is reducible. We will see shortly  that  $Der(\F_n, \L)$ 
contains the $P_n$-fixed line $B^1(\F_n, \L)$ and the $P_n$-invariant hyperplane 
$Der(\Ga_n, \L)$. The line does not intersect the hyperplane, nor it is 
complementary to it ($\L$ is not a field). We begin 
by describing the line. 

We have  seen that $Der(\F_n, \L)\cong Z^1(\F_n,\L)$. Consequently, 
$Der(\F_n, \L)$ contains $B^1(\F_n,\L)$, the Eilenberg-MacLane 
1-coboundaries. Since $C^0(\F_n, \L)\cong \L$ and $P_n$ acts trivially 
on $\L$, $P_n$ will also act trivially on $B^1(\F_n,\L)$. 

\begin{lem}
$B^1(\F_n,\L)$ is a free rank 1 submodule of $Z^1(\F_n, \L)$ with the basis 
$\sum_{i=1}^n (1-t_i) \frac{\D}{\D x_i}$. 
\end{lem}
\proof Recall that the coboundary $\del: C^0(\F_n, \L) \to 
 C^1(\F_n, \L)$ is given by
$$
\del \ell(x_i)= \ell - x_i \ell = \ell -  t_i \ell= (1-t_i)\ell
$$
But $(1-t_i)\ell= \ell \del 1(x_i)$, thus $\del$ is $\L$-linear 
and $B^1(\F_n,\L)= \L (\del 1)$. We conclude by observing that
$$
\del 1= \sum_{i=1}^n (1-t_i) \frac{\D}{\D x_i} \quad \qed
$$

We now describe the hyperplane. The element 
$x_{\infty} = x_1 \ldots x_n \in \F_n$ is fixed by $P_n$. We define 
$$
Der(\F_n, \L)^{\infty}:= \{ D \in Der(\F_n, \L) : D x_{\infty}= 0\}
$$   

\begin{lem}
(i) $Der(\F_n, \L)^{\infty}$ is a free summand of $Der(\F_n, \L)$ of rank $n-1$. 

(ii) The quotient map $\F_n \to \Ga_n$ induces an isomorphism 
$Der(\Ga_n,\L) \to Der(\F_n, \L)^\8$ of $P_n$-modules. 
 \end{lem}
\proof Let $\{y_1,..., y_n\}$ be the basis for $\F_n$ given by 
$y_i= x_1 \ldots x_i$, $1\le i\le n$. Then $Der(\F_n, \L)$ is free on $\frac{\D}{\D y_1},\ldots, 
\frac{\D}{\D y_n}$ and $Der(\F_n, \L)^{\infty}$ is free on 
$\frac{\D}{\D y_1},\ldots, \frac{\D}{\D y_{n-1}}$. The statement (ii) is clear. $\qed$ 

\begin{defn}
The reduced Gassner representation is the restriction of the action  of 
$P_n$ from $ Der(\F_n, \L)$ to $Der(\Ga_n,\L)$: 
$$
\rho: P_n \to Aut_{\L}( Der(\Ga_n,\L)).
$$
\end{defn}

\noindent We may represent $\rho(g)$, $g\in P_n$ as elements of 
$GL_{n-1}(\L)$ relative to the basis 
$\frac{\D}{\D y_1},\ldots, \frac{\D}{\D y_{n-1}}$. Observe that 
$B^1(\F_n, \L)$ does not intersect $Der(\Ga_n,\L)$, indeed
$$
\ell \del 1(x_{\infty})= \ell(1- t_1 \ldots t_n) \ne 0.
$$

\begin{rem}
We will see below that there exist homomorphism images of $Der(\F_n, \L)$ such that the image of 
$B^1(\F_n,\L)$ is contained in the image of 
$Der(\Ga_n,\L)$. Hence  $B^1(\F_n,\L)$ is not a complement to 
$Der(\Ga_n,\L)$. 
\end{rem}

Note also that there is a representation of $P_n$ on $\H^1(\F_n, \L)= Z^1(\F_n,\L)/B^1(\F_n,\L)$. 
We do not know whether or not $\H^1(\F_n,\L)$ is a  free $\L$-module. 

We now have 

\begin{defn}
Let $\al=(\al_1,..., \al_n)$ with $\al_j\in \C^*$, $1\le j\le n$ and $\M$ be an 
$\L$-module. Then the {\bf specialization} $\M_{\al}$ of $\M$ at $\al$ is defined by 
$\M_{\al}= \M \otimes_{\L} \C_{\al}$. Here $\C_{\al}$ is the complex line equipped 
with the $\L$-module structure $t_i z= \al_i z$, $z\in  \C$. 
\end{defn}

More concretely, $\M_{\al}$ is the quotient of $\M$ by the submodule of elements 
$\{(t_j -\al_j)m, 1\le j\le n, m\in \M\}$. 

Suppose that $T\in  End_{\L}(\M)$. Then $T$ induces an element $T_{\al}= T\otimes 1$ 
of $End(\M_{\al})$. Now assume that $\M$ is free on $m_1,...,m_n$. Then 
$m_1\otimes 1,..., m_n\otimes 1$ is a vector space basis for $\M_\al$. The matrix of 
$T_\al$ relative to this basis is obtained from a  matrix of $T$ relative to $m_1,...,m_n$ 
by substituting $\al_j$ for $t_j$, $1\le j\le  n$. 

Now we return to the case in hand. We have $\la_1,...,\la_n$ with 
$\la_1 + ... + \la_n=0$. Define $\al_j:= e^{2\pi i\la_j}$, $1\le j\le n$; 
whence $\al_1 ... \al_n=1$.

\begin{lem}
Suppose that $\al=(\al_1,...,\al_n)$ satisfies $\al_1 ... \al_n=1$. Then in the specialization  
$Der(\F_n,\L)_{\al}$ the image of the fixed line $B^1(\F_n,\L)$ is contained in the image of 
the invariant hyperplane $Der(\Ga_n,\L)$.  
\end{lem}
\proof $\del 1(x_{\infty})= 1- \al_1 ...\al_n =0$. $\qed$ 

\begin{cor}
The specialization  $Der(\Ga_n,\L)_{\al}$ contains a $P_n$-fixed line $B^1(\F_n,\L)_{\al}$. 
\end{cor}

Now we observe that $Z^1(\F_n,\L)_{\al}= Z^1(\F_n, \C_{\chi})$, the group of 
1-cocycles with values in the 1-dimensional module defined by $\chi(x_j)=\al_j, 
1\le j\le n$. Moreover
$$
Z^1(\Ga_n,\L)_{\al}= Z^1(\Ga_n, \C_{\chi}), 
$$
the group of $\C_{\chi}$-valued 1-cocycles  that annihilate $x_{\infty}$ and 
$$
B^1(\F_n,\L)_{\al}= B^1(\Ga_n, \C_{\chi}).
$$
We obtain

\begin{prop}
Suppose $\al=(\al_1,...,\al_n)$ satisfies $\al_1 ... \al_n=1$. Then the specialization of 
the reduced Gassner representation at $\al$ contains a $P_n$-invariant line. The quotient 
of the representation of $P_n$ by this line 
is  $\H^1(\Ga_n, \C_{\chi})$.  
\end{prop}

Theorem C follows.

\noindent Michael Kapovich: Department of Mathematics, University of 
Utah,  Salt Lake City, UT 84112, USA ; kapovich$@$math.utah.edu

\smallskip
\noindent John J. Millson: Department of Mathematics, University of 
Maryland, College Park, MD 20742, USA ; jjm$@$math.umd.edu

\end{document}